\documentclass{amsart}

\address{Michela Artebani, Departamento de Matem\'atica, Universidad de Concepci\'on, Casilla 160-C, Concepci\'on, Chile}

\email{martebani@udec.cl }

\address{Alessandra Sarti, Fachbereich f{\"u}r Mathematik, Johannes Gutenberg-Universit{\"a}t,
55099 Mainz, Germany}

\email{sarti@mathematik.uni-mainz.de}

\urladdr{http://www.mathematik.uni-mainz.de/$\sim$sarti}

\usepackage{amsmath }
\usepackage{bbm,dsfont}
\usepackage{mathrsfs}
\usepackage[all]{xy}
\usepackage{amsfonts}
\usepackage{ifthen}
\usepackage{rotating}
\usepackage{amssymb}
 
\usepackage{tikz}

\usepackage{enumerate}
 
\DeclareMathOperator{\rank}{rank}

\newcommand{\C}{\mathbb{C}}
\newcommand{\Z}{\mathbb{Z}}
\newcommand{\lra}{\longrightarrow}
\newcommand{\ra}{\rightarrow}

\newcommand{\PP}{\mathbb{P}}

 \newcommand{\Ps}{\mathbb P^}

\newtheorem{Lem}{Lemma}[section]
\newtheorem{theorem}{Theorem}[section]
\newtheorem{lemma}[theorem]{Lemma}
\newtheorem{pro}[theorem]{Proposition}
\newtheorem{cor}[theorem]{Corollary}

 \theoremstyle{definition}
\newtheorem{rem}[theorem]{Remark}
  \newtheorem{exa}[Lem]{Examples}
  \newtheorem{defin}[theorem]{Definition}
  
\begin{document}
\title{Non-symplectic automorphisms of order 3\\ on K3 surfaces}
\author{Michela Artebani \and Alessandra Sarti}
\date{\today}
\begin{abstract}
In this paper we study  K3 surfaces with a non-symplectic automorphism of order $3$. In particular, we classify the topological structure of the fixed locus of such automorphisms and we show that it determines the action on cohomology. This allows us to describe the structure of the moduli space and to show that it has three irreducible components.  
\end{abstract}



\maketitle
 
\pagestyle{myheadings} 
\markboth{Michela Artebani and Alessandra Sarti}{K3 surfaces with a non-symplectic automorphism of order three }
\setcounter{tocdepth}{1}

\section*{Introduction}
An automorphism on a K3 surface is called \emph{non-symplectic} if its natural representation on the vector space of holomorphic two-forms is not trivial.
In \cite{N1} and \cite{S} it was proved that a purely non-symplectic group of automorphisms is finite and cyclic.
All the orders of such groups have been determined in  \cite{MO} by Machida and Oguiso. In particular, the maximal order is known to be $66$ and $19$ if it is a prime.
Non-symplectic involutions have been studied in \cite{N3} and \cite{Z2}. In \cite{OZ2} Oguiso and Zhang showed that the moduli space of K3 surfaces with a non-symplectic automorphism of order $13, 17$ or $ 19$ is zero dimensional. The case of order $11$ has been studied by the same authors in \cite{OZ1}. 

The fixed locus of a non-symplectic involution is known to be either empty or the disjoint union of smooth curves and it is completely described.
Moreover, the action of the involution on the K3 lattice is well understood and only depends on the topology of the fixed locus.

In this paper we intend to give similar results for a non-symplectic automorphism of order $3$ on a K3 surface $X$ i.e.  $\sigma\in Aut(X)$ such that
$$\sigma^3=id \hspace{0.5cm} \mbox{and}\hspace{0.5cm} \sigma^*(\omega_X)=\zeta \omega_X,$$
where $\omega_X$ is a generator for $H^{2,0}(X)$ and $\zeta$ is a primitive $3$-rd root of unity. 

We prove that the fixed locus of $\sigma$ is not empty and it is the union of $n\leq 9$ points and $k\leq 6$ disjoint smooth curves, where all possible values of $(n,k)$ are given in Table \ref{fixtable}.
 
The automorphism $\sigma$ induces a non-trivial isometry $\sigma^*$ on the K3 lattice. The fixed sublattice for $\sigma^*$  is contained in the Picard lattice and its orthogonal complement is known to be a lattice over the Eisenstein integers. 
Here we prove that the isometry $\sigma^*$ only depends on the pair $(n,k)$.
The fixed lattice and its orthogonal complement have been computed for any $(n,k)$ and are listed in Table \ref{lat}.

The type $(n,k)$ gives a natural stratification of the moduli space of K3 surfaces with a non-symplectic automorphism of order $3$. As a consequence of the previous results, we prove that each stratum is birational to the quotient of a complex ball by the action of an arithmetic group. In particular, we prove that there are $3$ maximal strata of dimensions $9, 9, 6$, corresponding to $(n,k)=(0,1), (0,2), (3,0)$. The first two components also appear in \cite{K}, in particular the first one is known to be birational to the moduli space of curves of genus $4$.\\
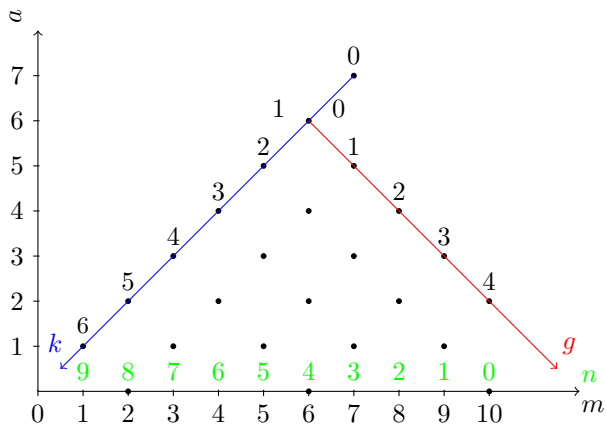
\begin{figure}
\begin{tikzpicture}[scale=.6]
\filldraw [black] 
(1,1) circle (1.5pt) node[below=-0.5cm]{\small{6}}
(2,2) circle (1.5pt) node[below=-0.5cm]{5}
(2,0) circle (1.5pt) 
(3,1) circle (1.5pt)
(3,3) circle (1.5pt) node[below=-0.5cm]{4}
(4,4) circle (1.5pt) node[below=-0.5cm]{3}
(4,2) circle (1.5pt)
(5,5) circle (1.5pt) node[below=-0.5cm]{2}
(5,3) circle (1.5pt)
(5,1) circle (1.5pt)
(6,6) circle (1.5pt) node[left=0.4cm, below=-0.4cm]{1} node[right=0.4cm,below=-0.4cm]{0}
(6,4) circle (1.5pt)
(6,2) circle (1.5pt)
(6,0) circle (1.5pt)
(7,7) circle (1.5pt)node[below=-0.5cm]{0}
(7,5) circle (1.5pt)node[ below=-0.5cm]{1}
(7,3) circle (1.5pt)
(7,1) circle (1.5pt)
(8,4) circle (1.5pt)node[ below=-0.5cm]{2}
(8,2) circle (1.5pt)
(9,3) circle (1.5pt)node[ below=-0.5cm]{3}
(9,1) circle (1.5pt)
(10,2) circle (1.5pt)node[ below=-0.5cm]{4}
(10,0) circle (1.5pt) node[ below=-0.5cm, green]{0}
 %(11,1) circle (0pt)node[ below=-0.5cm] {5}
(9,0) circle (0pt) node[ below=-0.5cm, green]{1}
(8,0) circle (0pt) node[ below=-0.5cm, green]{2}
(7,0) circle (0pt) node[ below=-0.5cm, green]{3}
(6,0) circle (0pt) node[ below=-0.5cm, green]{4}
(5,0) circle (0pt) node[ below=-0.5cm, green]{5}
(4,0) circle (0pt) node[ below=-0.5cm, green]{6}
(3,0) circle (0pt) node[ below=-0.5cm, green]{7}
(2,0) circle (0pt) node[ below=-0.5cm, green]{8}
(1,0) circle (0pt) node[ below=-0.5cm, green]{9}
; 
\draw plot[mark=*] file {data/ScatterPlotExampleData.data};
\draw[->] (0,0) -- coordinate (x axis mid) (12,0);
    \draw[->] (0,0) -- coordinate (y axis mid)(0,8);
    \foreach \x in {0,1,2,3,4,5,6,7,8,9,10}
        \draw [xshift=0cm](\x cm,1pt) -- (\x cm,-3pt)
         node[anchor=north] {$\x$};
          \foreach \y in {1,2,3,4,5,6,7}
        \draw (1pt,\y cm) -- (-3pt,\y cm) node[anchor=east] {$\y$};
    \node[below=0.2cm, right=3.5cm] at (x axis mid) {$m$};
    \node[below=-0.2cm, right=3.5cm, green] at (x axis mid) {$n$};
    \node[left=0.5cm, below=-2.6cm, rotate=90] at (y axis mid) {$a$};
 \draw[<-, blue](0.5,0.5)-- node[below=1.6cm,left=1.8cm]{$k$} (7,7);   
 \draw[<-, red](11.5,0.5)-- node[below=1.35cm,right=1.6cm]{$g$} (6,6);
  \end{tikzpicture} 
\caption{Fixed locus and fixed lattice}\label{fig}
\end{figure}
In the first section we recall some basic properties of non-symplectic  automorphisms of order three and their relation with lattices over the Eisenstein integers.

In section 2 we determine the structure of the fixed locus by generalizing a method in \cite{N3} and \cite{Ka}. More precisely, we determine algebraic relations between $n,k$ and two integers $m,a$ which identify the fixed lattice of $\sigma^*$.
All possible values of these two pairs of invariants are given in Table \ref{fixtable} and represented in Figure \ref{fig} (the analogous diagram for non-symplectic involutions can be found in \S 2.3, \cite{AN}).

In section 3 we prove that each such configuration occurs i.e. it can be realized by an order $3$ non-symplectic automorphism on a K3 surface. This is done by means of lattice theory, Torelli theorem and surjectivity theorem for the period map of K3 surfaces.
The list of the corresponding invariant lattices and of their orthogonal complements is given in Table \ref{lat}.

In section 4 we describe projective models for K3 surfaces with non-symplectic automorphisms of order $3$.
In particular, we show that for $k>1$ all of them have a jacobian elliptic fibration such that the automorphism acts as the identity on the basis.
Other models are given as complete intersections in $\Ps 4$, quartic surfaces and double sextics.

In the last section we describe the structure of the moduli space. 
The projective models given in section 4 show that the
moduli space is irreducible for given $n,k$ or, equivalently, that the action of the automorphism on cohomology is determined by $n,k$. 
This result and lattice theory allow us to identify the irreducible components of the moduli space. \\

\emph{Acknowledgements.} We would like to thank Igor Dolgachev, Shigeyuki Kond\=o, Bert van Geemen and Jin-Gen Yang for several stimulating discussions.

\section{Automorphisms and lattices}
Let $X$ be a K3 surface i.e. a simply connected smooth complex projective surface with a nowhere vanishing holomorphic two-form $\omega_X$. An automorphism $\sigma$ of $X$ is called \emph{non-symplectic} if its action on the vector space $H^{2,0}(X)=\C\omega_X$ is not trivial.
In this paper we are interested in non-symplectic automorphisms of order three i.e.
$$ \sigma^3=id \hspace{0.5cm}\mbox{ and }\hspace{0.5cm} \sigma^*(\omega_X)=\zeta \omega_X,$$
where $\zeta$ is a primitive $3$-rd root of unity.

The automorphism $\sigma$ induces an isometry $\sigma^*$ on $H^2(X,\Z)$ which preserves the Picard lattice $N_X$ and the transcendental lattice $T_X$ of $X$
$$N_X=\{x\in H^2(X,\Z): (x,\omega_X)=0\}\hspace{1cm}  T_X=N^{\perp}_X .$$  
We will denote by $N(\sigma)$ the invariant lattice
$$N(\sigma)=\{x\in H^2(X,\Z): \sigma^*(x)=x\}.$$

In order to describe the action of the automorphism on cohomology, we first need some preliminaries of lattice theory (see  \cite{N2})
Recall that the \emph{discriminant group} of a lattice $L$ is the finite abelian group $A_L=L^*/L$, where $L^*=Hom(L,\Z)$  and $A_L\cong A_{L^{\perp}}$.
\begin{defin} Let $\mathcal E=\Z[\zeta]$ be the ring of Eisenstein integers. A \emph{$\mathcal E$-lattice} is a pair $(L,\rho)$ where $L$ is an even lattice and $\rho$ an order three fixed point free isometry on $L$.
If $\rho$ acts identically on $A_L$ then $L$ will be called \emph{$\mathcal E^*$-lattice}.
\end{defin}
\begin{rem}\label{herm}
Any $\mathcal E$-lattice $L$ is clearly a module over $\mathcal E=\Z[\zeta]$ via the action
$$(a+\zeta b)\cdot x=ax+b\sigma^*(x),\hspace{0.7cm} a,b\in\Z, \ x\in L.$$ 
Since $\rho$ is fixed point free and $\mathcal E$ is a principal ideal domain, we have that $L$ is a free module over $\mathcal E$. Moreover, $L$ is equipped with a hermitian form $H:L\times L\ra \mathcal E$ such that $2 Re\,H(x,y)=(x,y)$. This is defined as
$$H(x,y)=\frac{1}{2}[(x,y)-\frac{\theta}{3}(x,\sigma^2(y)-\sigma(y))]$$
where $\theta=\zeta-\zeta^2$. 
\end{rem}
\begin{lemma}\label{inv}
\begin{enumerate}[$\bullet$]\ \\
\vspace{-0.3cm}
\item Any $\mathcal E$-lattice has even rank.
\item Any $\mathcal E^*$-lattice
is {\it 3-elementary} i.e. its discriminant group is isomorphic to $ (\Z/3\Z)^{\oplus a}$ for some $a$.
\end{enumerate}
\end{lemma}
\proof
The rank of a $\mathcal E$-lattice is equal to $2m$, where $m$ is its rank over $\mathcal E$.

Since $\rho$ is fixed point free, then $\rho^2+\rho+id=0$ on $L$.
If  $x\in A_{L}$  then $\rho(x)=x$ hence  $\rho^2(x)+\rho(x)+x=3x=0$. Thus $A_{L}$ is a direct sum of copies of $\Z/3\Z$.
\qed\\

Note that, according to Lemma \ref{inv}, to any $\mathcal E^*$-lattice $L$ we can associate the pair $m(L),a(L)$ where $2m(L)=\rank (L)$ and $a(L)$ is the minimal number of generators of $A_L$.\\

The following is a reformulation of \cite[Theorem 0.1]{N1} and \cite[Lemma 1.1]{MO} for order three automorphisms.
  
\begin{theorem}\label{ns}
Let $X$ be a K3 surface and $\sigma$ be a non-symplectic automorphism of order $3$ on $X$. Then:
\begin{enumerate}[$\bullet$]
\item  $N(\sigma)$ is a primitive $3$-elementary sublattice of $N_X$,
 \item $(N(\sigma)^{\perp},\sigma^*)$ is a $\mathcal E^*$-lattice
 \item $(T_X,\sigma^*)$ is a $\mathcal E$-sublattice of $N(\sigma)^{\perp}$.

\end{enumerate}
\end{theorem}

\begin{exa}\label{ex} In what follows we will adopt the standard notation for lattices: $U$ is the hyperbolic lattice and $A_n, D_n,  E_n$ are the negative definite lattices of rank $n$ associated to the corresponding root systems. The lattice $L(\alpha)$ is obtained multiplying by $\alpha$ the form on $L$.

\begin{enumerate}[$\bullet$]
 \item The lattices $U$, $E_8$, $A_2$, $E_6$,   $U(3)$ are $3$-elementary lattices with $a=0,0,1,1,2$ respectively.
\item If $L$ is a $3$-elementary lattice of rank $n$ with $a(L)=a$, then it can be proved that the scaled dual $L^*(3)$ is again a $3$-elementary lattice with $a(L^*(3))=n-a$. For example, the lattice $E_6^*(3)$ is $3$-elementary with $a=5$.
\item The lattice $A_2$ is a $\mathcal E^*$-lattice with order $3$ isometry : 
$$e\mapsto f,\ f\mapsto -e-f.$$
where $e^2=f^2=-2$, $(e,f)=1$.
 \item
The lattice $U\oplus U(3)$ is a $\mathcal E^*$-lattice with order three isometry :
$$e_1\mapsto e_1-f_1,\ e_2\mapsto -2e_2-f_2,$$
$$f_1\mapsto -2f_1+3e_1,\ f_2\mapsto f_2+3e_2,$$
where $e_i^2={f_i}^2=0$, $(e_1,e_2)=1$, $(f_1,f_2)=3$.
\item 
The lattice $U\oplus U$ has a structure of $\mathcal E^*$-lattice induced by the natural embedding $U\oplus U(3)\subset U\oplus U$.
\item The lattices $E_6$, $E_8$ are $\mathcal E^*$-lattices, see \cite[\S 2.6, Ch.2]{CS} for a description of two order three isometries on them.
\item The \emph{Coxeter-Todd lattice} $K_{12}$ is a  negative definite $\mathcal E^*$-lattice of rank $12$ with $a=6$
(see \cite{CS1}).   
In fact this is a \emph{unimodular lattice} over $\mathcal E$ i.e. $Hom(K_{12},\mathcal E)=K_{12}$. In \cite[Theorem 3]{F} W. Feit proved that this is the only unimodular $\mathcal E$-lattice of dimension $<12$ containing no vectors of norm $1$. 
 \end{enumerate}
\end{exa}

Hyperbolic $3$-elementary lattices have been classified by Nikulin \cite{N2} and  Rudakov-Shafarevich \cite{RS} (note that in this paper there is a misprint in the last condition of the theorem).
\begin{theorem}\label{RS}
An even hyperbolic $3$-elementary lattice $L$ of rank $r>2$ is uniquely determined by the integer $a=a(L)$. Moreover, given $r$ and $a\leq r$ such a lattice exists if and only if the following conditions are satisfied
$$
\begin{array}{ll}
r\equiv 0\ (mod\,  2)&\\
r\equiv 2\ (mod\, 4)& \mbox{ for } a\equiv 0\ (mod\, 2)\\
(-1)^{r/2-1}\equiv 3\ (mod\, 4)& \mbox{ for } a\equiv 1\ (mod\, 2)\\
r>a>0& \mbox{ for } r\not \equiv 2\ (mod\, 8)\\
\end{array}
$$
\end{theorem}
For $r=2$ binary forms are classified (see \cite{Bu}, \cite{CS}), in this case the unique definite even $3$-elementary lattices  are $A_2$ and $A_2(-1)$ (with $a=1$) and the only indefinite ones are $U$ and $U(3)$ (with $a=0$ and $a=2$).

\section{The fixed locus}
Let $\sigma$ be an order three non-symplectic automorphism of a K3 surface and $X^{\sigma}$ be its fixed locus. By \cite{Ca} the action of $\sigma$ can be locally linearized and diagonalized at $p\in X^{\sigma}$. Since $\sigma$ acts on $\omega_X$ as the multiplication by $\zeta$, there are two possible local actions for $\sigma$ 
$$ \left(\begin{array}{cc}
\zeta^2&0\\
0&\zeta^2
\end{array}
\right)~~~\mbox{or}~~~\left(\begin{array}{cc}
\zeta&0\\
0&1
\end{array}
\right)
.$$

In the first case $p$ is an isolated fixed point and $\sigma$ acts as the identity on $\PP T_p(X)$. 
In the second case $p$ belongs to a curve in the fixed locus (the line $x=0$). Hence
\begin{lemma}
The fixed locus of $\sigma$ is either empty or the disjoint union of $k$ smooth curves and $n$ points. 
  \end{lemma}

In this section we will relate the topological invariants $n, k$ with the lattice invariants $m,a$ defined in the previous section. Our approach generalizes a technique in \cite{Ka}. 
\begin{theorem}\label{fix}
The fixed locus of  $\sigma$ is not empty and it is the disjoint union of $n\leq 9$ points and $k\leq 6$ smooth curves with:
\begin{enumerate}[a)]

\item  one curve of genus $g\geq 0$ and $k-1$ rational curves or
\item $k=0$ and $n=3$.
\end{enumerate}
 Moreover, if $\rank N(\sigma)^{\perp}=2m$, then
$$m+n=10\hspace{0.5cm} \mbox{ and }\hspace{0.5cm} g=3+k-n\ \mbox{ (in\ case a))}.$$
  \end{theorem}
\proof
Let $\sigma$, $X$, $X^{\sigma}$ and $N(\sigma)$ as before. 
By the topological Lefschetz fixed point formula   
$$
\chi(X^{\sigma})=2+\mbox{Tr}(\sigma^*_{|H^2(X,\Z)})=2+\rank(N(\sigma))+m(\zeta+\zeta^2)
.$$
Since $\zeta+\zeta^2+1=0$ and $\rank(N)=22-2m$ this gives 
\begin{equation}\label{eq1}
\chi(X^{\sigma})=\sum_{i=1}^k \mathcal X(C_i)+n= 3(8-m)
,\end{equation}
where $C_i$ are the smooth curves in the fixed locus. 
 Besides, the holomorphic Lefschetz formula (Theorem 4.6, \cite{AS}) gives the equality
\begin{equation}\label{eq2}
 1+\mbox{Tr}(\sigma^*_{|H^{2,0}(X) })= \frac{1}{6\zeta}(\sum_{i=1}^k\mathcal X(C_i)-2n)\Longleftrightarrow 
2n-\sum_{i=1}^k \mathcal X(C_i)=6.
\end{equation}
The equations \eqref{eq1} and \eqref{eq2} imply $m+n=10$, hence also $n\leq 9$. Case $b)$ immediately follows from  (\ref{eq2}).

By the Hodge index theorem the Picard lattice $N_X$ is hyperbolic. This implies that the fixed locus contains at most one curve $C$ of genus $g>1$ and that in this case the other curves in the fixed locus are rational (since they belong to the orthogonal complement of $C$).  

If the fixed locus of $\sigma$ contains at least two elliptic curves $C_1, C_2$, then they are linearly equivalent and $|C_1|:X\lra\Ps 1$ gives an elliptic fibration on $X$. The induced action of $\sigma$ on $\Ps 1$ is not trivial since otherwise $\sigma$ would act as the identity on the tangent space at a point in $C_1$. Hence $\sigma$ has exactly $2$ fixed points in $\Ps 1$, corresponding to the fibers $C_1,C_2$ and there are no other fixed curves or points in the fixed locus.
This contradicts equality (\ref{eq2}), hence the fixed locus contains at most one elliptic curve. This completes the proof of $a)$.

As a consequence, if $g$ is the biggest genus of a curve in the fixed locus, (\ref{eq2}) can be written as $g=3+k-n$.

The inequality $k\leq 6$ follows from (\ref{eq2}) and the Smith inequality (see \cite{Br})
$$\sum_{i\geq 0} \dim H^i(X^{\sigma},\Z_3)\leq \sum_{i\geq 0} \dim H^i(X,\Z_3).$$\qed  
\begin{cor}
An order $3$ automorphism on a K3 surface is symplectic if and only if its fixed locus is given by $6$ points.
\end{cor}
\proof
An order $3$ symplectic automorphism on a K3 surface has exactly 6 fixed points by \cite[\S 5]{N1}. The converse follows from Theorem \ref{fix}.
\qed\\

Let $a=a(N(\sigma))$, then a more refined application of Smith exact sequences gives

\begin{pro}\label{a}
$\dim H_{*}(X)-\dim H_{*}(X^{\sigma})=2a+m.$
\end{pro}
\proof In what follows we will write $\sigma$ for $\sigma^*$.
Let $g=id+\sigma+\sigma^{2}$ and $h=id-\sigma$ on $H^2(X,\Z)$.
If $L_+=ker(h)$ and $L_-=ker(g)$, then
$$3^a=o(H^2(X,\Z)/L_+\oplus L_-).$$
Since $h^2=id-2\sigma+\sigma^2=g$ over $\Z_3$, then the image of $L_+\oplus L_-$ under the homomorphism
$$c:H^2(X,\Z)\lra H^2(X,\Z_3)$$
coincides with $c(L_-)$. Hence
$$a=\dim H^2(X,\Z_3)-\dim c(L_-).$$
We now express this number in a different way by applying Smith exact sequences over $\Z_3$ (see Ch.III, \cite{Br}).
In the rest of the proof the coefficients will be in $\Z_3$.

Let $C(X)$ be the chain complex of $X$ with coefficients in $\Z_3$. 
%$\mathcal X(\cdot)$ be the Euler characteristic of a complex and $\mathcal X(\cdot,\cdot)$ be the Euler characteristic of a pair. 
%We will identify $X^{\sigma}$ with its image in $S=X/\sigma$.
The automorphisms $g,h$ act on $C(X)$ and give chain subcomplexes $g C(X)$ and $h C(X)$. We denote by $H_i^{g}(X),\ H_i^{h}(X)$ the associated homology groups with coefficients in $\Z_3$ as in \cite[Definition 3.2]{Br} and by $\mathcal X^g(X), \mathcal X^h(X)$ the corresponding Euler characteristics. By \cite[Proposition 3.4]{Br} we have
$$H_i^g(X)\cong H_i(S,X^{\sigma}), $$
where the second term is the homology of the pair $(S,X^{\sigma})$ where $S=X/\sigma$ and $X^{\sigma}$ is identified to its image in $S$.

Let $\rho=h^i$ and $\bar \rho=h^{3-i},$ with $i=1,2$. Then we have the exact triangles (\cite[Theorem 3.3 and (3.8)]{Br})
$$\xymatrix{
(3)\label{t1}\ \ & &H(X) \ar[dl]_{\rho_*} &   \\
&H^{\rho}(X)\ar[rr]_{}& &H^{\bar\rho}(X)\oplus H(X^{\sigma})\ar[ul]_{i_*}\\  
(4)\label{t2}\ \ & &H^{h}(X) \ar[dl]_{h_*} &   \\
&H^{g}(X)\ar[rr]_{}& &H^{g}(X) \ar[ul]_{i_*}} 
$$
where $h_*, i_*, \rho_*$ have degree $0$ and the horizontal arrows have degree $-1$.
The two triangles \eqref{t1} and (4) are called Smith sequences, in particular they induce two exact sequences 
$$0 \ra H_3^{g}(X)\stackrel{\gamma_3}{\ra} H_2^{h}(X)\oplus H_2(X^{\sigma})\stackrel{\alpha_2}{\ra} H_2(X)\stackrel{\beta_2}{\ra} H^{g}_2(X)\stackrel{\gamma_2}{\ra} H^{h}_1(X)\oplus H_1(X^{\sigma})\ra 0 $$
$$0\ra H_3^{h}(X)\stackrel{\gamma'_3}{\ra} H_2^{g}(X)\oplus H_2(X^{\sigma})\stackrel{\alpha'_2}{\ra} H_2(X)\stackrel{\beta'_2}{\ra} H^{h}_2(X)\stackrel{\gamma'_2}{\ra} H^{g}_1(X)\oplus H_1(X^{\sigma})\ra 0,$$

From exact sequence \eqref{t1} with $\rho=\sigma$ and sequence (4)   follow (\cite[Theorem 4.3]{Br}) the equalities of Euler characteristics
$$\mathcal X(X)-\mathcal X(X^{\sigma})=\mathcal X^g(X)+\mathcal X^h(X)=3\mathcal X^g(X).$$
%where the maps $\alpha_2,\alpha'_2$ are induced by the sum of inclusions
%$$i:h^iC(X)\oplus C(X^{\sigma})\lra C(X),\ i=1,2 $$
%and the maps $\beta_2, \beta'_2$ are induced by the maps
%$$h^i: C(X)\lra h^iC(X),\ i=1,2. $$
Then, from the exactness of Smith sequences, Lemma \ref{chi} and Lemma \ref{im} below we have
$$\begin{array}{ccl}\dim H_{*}(X)-\dim H_{*}(X^{\sigma})&=&\mathcal X^g(X)+\mathcal X^h(X)-2\dim H_1(X^{\sigma})\\
 &=&\dim Im(\beta_2)+\dim Im(\beta'_2)\\
 &=&2a+ m.\end{array}$$

\qed
\begin{lemma}\label{chi}
$\mathcal X^g(X)+\mathcal X^h(X)=\sum_{i=1}^2 (-1)^i (\dim H^g_i(X)+\dim H^h_i(X)).$
 \end{lemma}
\proof
From the exact sequence for the pair $(S,X^{\sigma})$ and sequences $(3)$, $(4)$ it follows
$$\dim H^g_0(X)=\dim H^h_0(X)=0,$$
$$ \dim H^g_3(X)=\dim H^h_3(X)=\dim H^g_4(X)=\dim H^h_4(X).$$
This immediately implies the statement.
\qed
\begin{lemma}\label{im}
 $Im(\alpha_2)=c(L_-)$, $\dim Im(\alpha_2)-\dim Im(\alpha'_2)=m$.
\end{lemma}
\proof
By definition of the Smith exact sequence $Im(\alpha_2)\subset c(L_-)$.
Conversely, if $x\in c(L_-)$ then $\alpha'_2(\beta_2(x)\oplus 0)=g(x)=0$, hence $(\beta_2(x)\oplus 0)\in Im(\alpha'_2)$.
By definition the projection of $\gamma_3$ on the second factor is the boundary homomorphism of the sequence of the pair $(S,X^{\sigma})$ and this is injective since $H_3(S)=0$. It follows that $\beta_2(x)=0$ i.e. $x\in ker(\beta_2)=Im(\alpha_2)$.

By the two exact sequences above and the homology vanishing in Lemma \ref{chi} it follows that
$$\dim Im(\alpha_2)-\dim Im(\alpha'_2)=h_2^h(X)-h_2^g(X)$$
(in fact $\dim H_1^g(X)=\dim H_1^h(X)$). Because of sequence $(4)$ and Proposition \ref{fix} this also equals
$\mathcal X^g(X)=m$. 
\qed
\begin{cor}\label{amg} If $k=0$ then $m=a$, otherwise $2g=m-a$. 
 \end{cor}
\proof
It follows from equation (\ref{eq1}) in the proof of Proposition \ref{fix} and from Proposition \ref{a} since
$$\dim H_{*}(X)=\mathcal X(X)=24,$$
$$\dim H_{*}(X^{\sigma})-\mathcal X(X^{\sigma})=2h^1(X^{\sigma})=4g.$$
\qed\\
 Note that $g$, $m$ and $a$ are functions of $n,k$:
$$\begin{array}{l}
g(n,k)=3+k-n\\
m(n,k)=10-n\\
a(n,k)=n+4-2k.
\end{array}$$

\begin{theorem}\label{fixt}
  The fixed locus of $\sigma$ contains $n$ points and $k$ curves where $n,k$ are in the same row of Table \ref{fixtable}.
\begin{table}
 \begin{tabular}{c|c|c|c|c}
$n$& $k$& $g(n,k)$ &     $m(n,k)$ & $a(n,k)$\\
\hline
$0$ &$1,2$& $4,5$ &   $10$  & $2,0 $\\
\hline
$1$ & $1,2$& $3,4$  &   $9$ & $3,1 $\\
\hline
$2$& $1,2$& $2,3$&    $8$ & $4,2$\\
\hline
$3$& $0,1,2,3$& $  \emptyset,1,2,3$&    $7$ & $7,5,3,1$\\
\hline
$4$& $1,2,3,4$& $0,1,2,3$&   $6$ & $6,4,2,0$\\
\hline
$5$& $2,3,4$& $0,1,2 $&    $5$ & $5,3,1$\\
\hline
$6$& $3,4$& $0,1$ &    $4$ & $4,2 $\\
\hline
$7$& $4,5$& $0,1$  &   $3$ & $3,1$\\
\hline
$8$ & $5,6$& $0,1$  &    $2$ & $2,0$\\
\hline
$9$ & $6$& $0$ &   $1$ & $1$
\end{tabular}\\
\ \\
\ \\
\caption{Fixed locus and fixed lattice}\label{fixtable}
\end{table}
\end{theorem}
 \proof
By Theorem \ref{RS} and Corollary \ref{amg} we get
$$a\leq min (m, 22-2m),\ \ a=0 \mbox { or }a=22-2m \Longrightarrow m\equiv 2\ (mod\, 4).$$
Then the result follows from Theorem \ref{fix} and Corollary \ref{amg}.
\qed

\section{Existence}

Let $(T,\rho)$ be an $\mathcal E^*$-lattice of signature $(2,n-2).$ Assume  that $T$ has a primitive embedding in the K3 lattice $L_{K3}$ and let $N$ be its orthogonal complement in $L_{K3}$.
 %We  consider the eigenspace 
%$$V=\{z\in T\otimes \C: \rho(z)=\zeta z\}$$ and the period domain 
%for K3 surfaces with this action of $\sigma^*$ is given by
In $L_{K3}\otimes \C$ consider the period domain
$$B_{\rho}=\{\omega\in \mathbb P(T\otimes \C): (\omega,\bar \omega)>0,\rho(\omega)=\zeta \omega\}.$$
If $\omega\in B_{\rho}$ is a generic point, then by the surjectivity of the period map (\cite{Ku}, \cite{PP}) there exists a K3 surface $X=X_{\omega}$ with a marking $\phi: H^2(X,\Z)\rightarrow L_{K3}$ such that $\mathbb P(\phi_{\C}(\omega_X))=\omega$.

\begin{pro}\label{existence}
\begin{table}
\begin{tabular}{c|c|c|c}
$n$ & $k$&  $T(n,k)$&$N(n,k)$\\
\hline
$0$ &$1$&  $U\oplus U(3)\oplus E_8\oplus E_8$& $U(3)$\\ 
& $2$&  $U\oplus U\oplus E_8\oplus E_8$& U\\
\hline
$1$ &$1$& $U\oplus U(3)\oplus E_6\oplus E_8$& $U(3)\oplus A_2$\\ 
 & $2$& $U\oplus U\oplus E_6\oplus E_8$& $U\oplus A_2$\\
\hline 
$2$ &$1$&  $U\oplus U(3)\oplus E_6\oplus E_6$ &$U(3)\oplus A_2^2$\\
 &$2$& $U\oplus U\oplus E_6\oplus E_6$&$U\oplus A_2^2$\\
\hline 
$3$ &$0$& $U\oplus U(3)\oplus A_2^5$& $U(3)\oplus E_6^*(3)$\\ 
 & $1$&  $U\oplus U\oplus A_2^5$ & $U(3)\oplus A_2^3$\\ 
 &$ 2$& $U\oplus U(3)\oplus A_2\oplus E_8$& $U\oplus A_2^3$\\ 
 &$3$& $U\oplus U\oplus A_2\oplus E_8$& $U\oplus E_6$\\  
\hline
$4$ &$1$&  $U\oplus U(3)\oplus A_2^4$& $U(3)\oplus A_2^4$\\ 
 &$2$& $U\oplus U\oplus A_2^4$&$U\oplus A_2^4$\\
 &$3$&  $U\oplus U(3)\oplus E_8$& $U\oplus E_6\oplus A_2$\\ 
 &$4$&  $U\oplus U\oplus E_8$&$U\oplus E_8$\\
\hline
$5$ &$2$&  $U\oplus U(3)\oplus A_2^3$&$U\oplus A_2^5$\\ 
 &$3$& $U\oplus U(3)\oplus E_6$&$U\oplus A_2^2\oplus E_6$\\ 
&$4$& $U\oplus U\oplus E_6$&$U\oplus E_8\oplus A_2$\\
\hline
$6$ &$3$& $U\oplus U(3)\oplus A_2^2$&$U\oplus E_6\oplus A_2^3$\\ 
&$4$& $U\oplus U\oplus A_2^2$&$U\oplus E_6^2$\\ 
\hline
$7$ &$4$& $U\oplus U(3)\oplus A_2$&$U\oplus E_6\oplus E_6\oplus A_2$\\ 
&$5$& $U\oplus U\oplus A_2$&$U\oplus E_6\oplus E_8$\\ 
\hline
$8$ &$5$&  $U\oplus U(3)$&$U\oplus E_6\oplus E_8\oplus A_2$\\ 
&$6$&  $U\oplus U$&$U\oplus E_8\oplus E_8$\\ 
\hline
$9$ &$6$&  $A_2(-1)$&$U\oplus E_8\oplus E_8\oplus A_2$\\ 
\end{tabular}\\
\ \\
\ \\
\caption{The lattices $T(n,k),\ N(n,k)$.}\label{lat}
\end{table}
The K3 surface $X$ admits an order three non-symplectic automorphism $\sigma$ such that $\sigma^*=\rho$ on $T$ up to conjugacy and
$$N(\sigma)=N_X\cong N,\ N(\sigma)^{\perp}=T_X\cong T.$$
 
 \end{pro}
\proof
We consider the isometry on $N\oplus T$ defined by $(x,y)\mapsto(x,\rho(y))$. Since $\rho$ acts as the identity on $A_T$, then this isometry  defines an isometry $\bar\rho$ on $L_{K3}$.
If $\omega\in B_{\rho}$, then $\bar\rho(\omega)=\zeta \omega$ and, if $\omega$ is generic, we can assume that there are no roots in $T\cap \omega^{\perp}$. Then the statement follows from  \cite[Theorem 3.10]{Na}.
\qed

\begin{pro}\label{transc}
For any $n,k$ in the same row of Table \ref{fixtable} there exists a unique $3$-elementary lattice $T(n,k)$ of  signature $(2,2m(n,k)-2)$ with $a=a(n,k)$.
This lattice is a $\mathcal E^*$-lattice and has a unique primitive embedding in $L_{K3}$.\\
The lattice $T(n,k)$ and its
orthogonal complement $N(n,k)$ in $L_{K3}$ are given in Table \ref{lat}.  
\end{pro}
\proof
It can be easily checked that any $m,a$ in Table \ref{fixtable} is realized by one of the lattices in Table \ref{lat}.  
Any lattice $T$ in Table \ref{lat} is a direct sum of $\mathcal E^*$-lattices in Examples \ref{ex}, hence it is also a  $\mathcal E^*$-lattice (by taking the direct sum of the isometries on each factor).

By \cite[Theorem 1.12.2]{N2} there exists a primitive embedding of $T$ in $L_{K3}$ if and only if there exists a hyperbolic $3$-elementary lattice $N$ of rank $22-2m(n,k)$ and $a=a(n,k)$.
This follows from Theorem \ref{RS}.

Moreover, by \cite[Corollary 1.12.3]{N2}, the embedding of $N$ in $L_{K3}$ is unique, hence the embedding of $T$ in $L_{K3}$ is also unique.
\qed

\begin{theorem}\label{clas}
For any $n,k$ in the same row of Table \ref{fixtable} there exists a K3 surface $X$ with a non-symplectic automorphism $\sigma$ of order three such that 
\begin{enumerate}[$\bullet$]
\item $X^{\sigma}$ contains $n$ points, $k$ curves and the biggest genus of a fixed curve is $g(n,k)$,
\item $N_X=N(\sigma)\cong N(n,k)$ and  $T_X\cong T(n,k)$.
\end{enumerate}
\end{theorem}  
\proof
It follows from Proposition \ref{existence}, Proposition \ref{transc} and Theorem \ref{fix}.
\qed\\

\noindent In what follows we will call $X_{n,k}$ and $\sigma_{n,k}$ a K3 surface and an automorphism as in Theorem \ref{clas}.
\section{Examples and projective models}
We will start studying elliptic fibrations (see \cite{M} for basic definitions and properties) on K3 surfaces with a non-symplectic automorphism  of order $3$ .
The possible Kodaira types for stable fibers and the action of the automorphism on them is described by the following result in \cite{Z3}.
\begin{lemma}\label{ell}
Let $X$ be a K3 surface with a non-symplectic automorphism $\sigma$ of order three and $f:X\ra \Ps 1$ be an elliptic fibration. If $F$ is a singular fiber of $f$ containing at least one $\sigma$-fixed curve,
then it is one of the following Kodaira types:

\begin{enumerate}[$\bullet$]
\item  \textbf{$IV$:} $F=F_1+F_2+F_3$ and $F_1$ is the only fixed curve in $F$.
\item \textbf{$I_n$} with \textbf{$n=3,6,9,12,15,18$:} $F=F_1+\cdots +F_n$  where $F_i\cdot F_{i+1}=F_n\cdot F_1=1$ and the $\sigma$-fixed curves in $F$ are $F_1,F_4,\dots,F_{n-2}$.
\item \textbf{$I^*_{n-5}$} with \textbf{$n=5,8,11,14,17$:} 
$$F=F_1+F_2+2(F_3+\cdots + F_{n-2})+F_{n-1}+F_n$$
where $F_1\cdot F_3=F_i\cdot F_{i+1}=F_{n-2}\cdot F_n=1$, $2\leq i\leq n-2$ and the $\sigma$-fixed curves are $F_3, F_6,F_9, \dots, F_{n-2}$.
\item \textbf{$IV^*$:} 
$F=3F_1+2F_2+F_3+2F_4+F_5+2F_6+F_7$
and $F_1$ is the only $\sigma$-fixed curve in $F$.
\item \textbf{$III^*$:}
$F=4F_1+2F_2+3F_3+2F_4+F_5+3F_6+2F_7+F_8$
and $F_1, F_5, F_8$ are the only $\sigma$-fixed curves in $F$.
\item \textbf{$II^*$:}
$F=6F_1+3F_2+4F_3+2F_4+5F_5+4F_6+3F_7+2F_8+F_9$
where $F_1, F_7$ are the only $\sigma$-fixed curves in $F$.
\end{enumerate}
\end{lemma}
We now show that any $X_{n,k}$ admits an invariant elliptic fibration if the fixed locus contains more than two curves  
\begin{pro}\label{k>1}
If $k>1$ then the K3 surface $X_{n,k}$ is isomorphic to a jacobian elliptic fibration with Weierstrass model
$$y^2=x^3+ p_{12}(t)$$
where $p_{12}(t)$ has exactly the following multiple roots
\begin{enumerate}[$\bullet$]
\item $n$ double roots if $k=2$.
\item one $4$-uple root and $n-3$ double roots if $k=3$
\item one $5$-uple root and $n-4$ double roots if $k=4$
\item one $5$-uple root, one $4$-uple root and $n-7$ double roots for $k=5$
\item two $5$-uple roots and $n-8$ double roots for $k=6$.
\end{enumerate}
In these coordinates the non-symplectic automorphism $\sigma_{n,k}$  is 
$$(x,y,t)\mapsto (\zeta x,y,t).$$
Conversely, for any $k>1$, a jacobian fibration with the above properties is a K3 surface and the automorphism $\sigma_{n,k}$ has fixed locus of type $(n,k)$. 
\end{pro}
\proof
If $k=2$ then by Proposition \ref{transc} the Picard lattice of $X=X_{n,k}$ is isomorphic to $U\oplus A_2^{n}$. Hence $X$ has a jacobian elliptic fibration $f:X\ra \Ps 1$ with $n$ reducible fibers with dual graph $\tilde A_2$.

The isometry $\sigma^*$ acts as the identity on $N_X$, thus it preserves $f$ and any of its reducible fibers. Since either there are $n>2$ reducible fibers or $g(n,k)>1$  then $\sigma$ acts as the identity on $\Ps 1$ i.e. any fiber of $f$ is $\sigma$-invariant and the section is in the fixed locus of $\sigma$. 
In particular, any fiber of $f$ has an automorphism of order $3$ with a fixed point.

The Weierstrass model of $f$ can be given by an equation of type
$$y^2=x^3+a(t)x+b(t).$$  
By the previous remark the functional invariant $j(t)=a(t)^3/\Delta(t)$ is equal to zero i.e. $a(t)\equiv 0$ and the action of $\sigma$ is $(x,y,t)\mapsto (\zeta x,y,t)$.
Moreover, the reducible fibers of $f$ are of Kodaira type $IV$ and there are $12-2n$ other singular fibers of type $II$.
Up to a change of coordinates in $\Ps 1$, we can assume that all singular fibers are in $\Ps 1\backslash\infty$ i.e. that $\deg(b(t))=12$. Note that simple roots of $b(t)$ give type $II$ fibers and double roots give type $IV$ fibers (see \cite{M} or \cite{Ko}). This gives the result.

For $k>2$ we can still see from $N(n,k)$ in Table \ref{fixtable} that $X$ has a jacobian elliptic fibration with fibers of type $\tilde A_2$, $IV^*$ and $II^*$.
By the previous arguments and Lemma \ref{ell} we find that their j-invariants are zero and we can write the Weierstrass equations as before, recalling that $4$-uple roots of $b(t)$ give type $IV^*$ fibers and $5$-uple roots give type $II^*$ fibers. 
 
 The last statement follows easily by Lemma \ref{ell}.
 \qed\\
 
 \noindent If $g(n,k)>0$, we will denote by $C_{n,k}$ the  curve  in the fixed locus of $\sigma_{n,k}$ with such genus.
\begin{cor}\label{hyp}
 If $k>1$ then $C_{n,k}$ is hyperelliptic.
 \end{cor}
 \proof
If $k>1$ then, by Proposition \ref{k>1}, $X_{n,k}$ has a jacobian elliptic fibration and the section is in the fixed locus. Since $\sigma_{n,k}$ fixes $3$ points on each fiber of the elliptic fibration, then $C$ is a double section and $f_{|C}$ is a double cover to $\Ps 1$. \qed\\

\noindent Let $\phi_{n,k}: X_{n,k}\lra \Ps {g(n,k)} \cong \mathbb |C_{n,k}|^{\vee} $ be the morphism associated to $|C_{n,k}|$.
 \begin{lemma}\label{proj}
There is a projective transformation $\tilde \sigma_{n,k}$ of $\Ps g$ which preserves $Im(\phi_{n,k})$ and such that $\phi_{n,k} \circ \sigma_{n,k}=\tilde \sigma_{n,k} \circ \phi_{n,k}$. For a suitable choice of coordinates  
$$\tilde \sigma_{n,k}(x_0,\dots,x_{g-1},x_g)=(x_0,\dots,x_{g-1},\zeta x_g).$$ 

\end{lemma}
\proof
Since $\sigma_{n,k}$ fixes $C_{n,k}$, then it preserves $|C_{n,k}|$, hence it induces a projectivity $\tilde \sigma_{n,k}$ of $|C_{n,k}|^{\vee}$  which fixes pointwisely the hyperplane $H$ such that $\phi_{n,k}^{-1}(H)=C_{n,k}$.
If we choose coordinates such that $H=\{x_g=0\}$, then $\tilde \sigma_{n,k}$ is of the above form.
\qed

\begin{rem}
 If $g(n,k)=1$ then $\phi_{n,k}$ is an elliptic fibration. The automorphism $\tilde\sigma_{n,k}$ has exactly two fixed points $p,q\in \Ps 1$ such that $\phi_h^{-1}(p)=C_{n,k}$ and $\phi_h^{-1}(q)$ is a reducible fiber. By \cite{D} this is of Kodaira type $I_0$  if $n=3$  and $I^*_{3(n-4)}$ for $n=4,\dots,8.$
 \end{rem}

 \begin{rem}
 The curve $C_{0,2}$ is a hyperelliptic curve of genus $5$ and by \cite{SD} $\phi_{0,2}:X_{0,2}\ra \Ps 5$ is a degree two morphism onto a cone over a rational normal quartic branched along a cubic section $B$ and the vertex of the cone. Since the branch curve is invariant for the action of $\tilde \sigma$, then $B$ has an equation of the form
$$F_3(x_0,\dots,x_4)+bx_5^3$$
where $b$ is non zero and $(0,\dots,0,1)$ is the vertex of the cone.
Conversely, the double cover of the cone branched along the generic section with this equation and the vertex is a K3 surface with an automorphism $\sigma_{0,2}$.

The surface $X_{0,2}$ can be also obtained as the triple cover of the Hirzebruch surface $\mathbb F_6$ branched along the disjoint union of the exceptional curve $e$ with $e^2=-6$ and a curve in $\mid 12f+2e \mid$, where $f$ is the class of a fiber. 
 \end{rem}

In what follows $F_i, G_i$ will denote homogeneous polynomials of degree $i$ and $b,c,d$ are non-zero complex numbers.

\begin{pro}\label{ic}
If $(n,k)=(0,1), (3,0) , (3,1)$ then $X_{n,k}$ is isomorphic to the complete intersection of a quadric and a cubic in $\Ps 4$ with equations of the form
 
$$
\begin{array}{l}
\bullet \hspace{0.7cm} X_{0,1}:\  \left\{
\begin{array}{ll}
F_2(x_0,\dots,x_3)=0\\
F_3(x_0,\dots,x_3)+bx_4^3=0 
\end{array}\right. \\
 \\
\hspace{0.9cm}   \sigma_{0,1}(x_0,\dots,x_3,x_4)=(x_0,\dots,x_3,\zeta x_4) \\
 
\ \\
\bullet \hspace{0.7cm} X_{3,0}:\ 
\left\{
\begin{array}{ll}
F_2(x_0,x_1)+b x_2 x_3+c x_2 x_4=0 \\
F_3(x_0,x_1)+d x_2^3+G_3(x_3,x_4)+x_2F_1(x_0,x_1)G_1(x_3,x_4)=0
\end{array}\right. \\
\ \\
\hspace{0.9cm} \sigma_{3,0}(x_0,x_1,x_2,x_3,x_4)=(x_0,x_1,\zeta^2x_2,\zeta x_3, \zeta x_4).\\
\ \\
\bullet \hspace{0.7cm} X_{3,1}:\ \left\{\begin{array}{ll}
x_3F_1(x_0,x_1,x_2)+x_4G_1(x_0,x_1,x_2)=0\\
F_3(x_0,x_1,x_2)+G_3(x_3,x_4)=0
\end{array}\right.\\
\ \\
\hspace{0.9cm} \sigma_{3,1}(x_0,x_1,x_2,x_3,x_4)=(x_0,x_1,x_2,\zeta x_3,\zeta x_4).
\end{array}$$
Conversely, for generic $F_i, G_j, b, c, d$, the above equations define K3 surfaces and the fixed locus of $\sigma_{n,k}$ is of type $(n,k)$.
 \end{pro}
   \proof
 The fixed locus of $\sigma_{0,1}$ only contains the curve $C_{0,1}$ of genus $4$.
Assume that $E$ is an elliptic curve on $X$ intersecting $C_{0,1}$, then $\sigma_{0,1}$ preserves $E$ and has exactly $3$ fixed points on it by Riemann Hurwitz formula, thus $(C_{0,1},E)=3$.
By \cite[Theorem 5.2]{SD} this implies that the morphism $\phi_{0,1}$ is an embedding of $X_{0,1}$ in $\Ps 4$ as the complete intersection of a quadric and a cubic hypersurface. Since both hypersurfaces are $\tilde \sigma_{0,1}$-invariant, then they have an equation of the above form and $\sigma_{0,1}=\tilde \sigma_{0,1}$.

According to Proposition \ref{transc}, the Picard lattice of $X_{3,0}$ is isomorphic to $U(3)\oplus E_6^*(3)$. Let $e,f$ be the standard basis of $U(3)$ and $h=e+f$. The morphism associated to $h$ is an embedding as the intersection of a quadric and a cubic hypersurface in $\Ps 4$.
Moreover, $\sigma_{3,0}$ induces a projectivity $\tilde \sigma$ of $\Ps 4$. Since $\sigma_{3,0}$ has only isolated fixed points, then $\tilde\sigma$ has no fixed hyperplane or planes. Hence, up to a choice of coordinates, we can assume that $\sigma_{3,0}$ and the two hypersurfaces have equations of the above form.

According to Proposition \ref{transc}, the Picard lattice of $X_{3,1}$ is isomorphic to $U(3)\oplus A_2^3$. Let $e,f$ be the basis of $U(3)$ as before and $e_i, e'_i$, $i=1,2,3$ be the standard basis of $A_2^3$.  
The class $h=2e+f-\sum_{i=1}^3 (e_i+e'_i)$ gives an embedding as the complete intersection of a quadric and a cubic hypersurface.
The equations of $X_{3,1}$ and $\sigma_{3,1}$ can be determined as before.

The last statement follows from easy computations.
\qed

\begin{rem}\label{g4}
It is clear from the equations in Proposition \ref{ic} that the surface $X_{0,1}$ is the triple cover of a quadric branched along a curve of genus $4$. Also the converse is true, as proved by Kond\=o in \cite{K}. In fact \cite[Theorem 1]{K} states that the moduli space of K3 surfaces $X_{0,1}$ is birational to the coarse moduli space of curves of genus $4$.

\end{rem}

\begin{pro}\label{qu}
If $(n,k)=(1,1), (4,1)$ then  $X_{n,k}$ is isomorphic to a smooth quartic in $\Ps 3$
with equations of the form
$$
\begin{array}{l}
\bullet \hspace{0.7cm} X_{1,1}:\   
F_4(x_0,x_1,x_2)+F_1(x_0,x_1,x_2)x_3^3=0\\
\ \\
\hspace{0.9cm}   \sigma_{1,1}(x_0,\dots,x_3)=(x_0,x_1,x_2,\zeta x_3) \\
\ \\
\bullet \hspace{0.7cm} X_{4,1}:\ 
F_4(x_0,x_1)+F_3(x_2,x_3)F_1(x_0,x_1)=0 \\
\ \\
\hspace{0.9cm} \sigma_{4,1}(x_0,x_1,x_2,x_3)=(x_0,x_1,\zeta x_2,\zeta x_3).\\

\end{array}$$
Conversely, for generic $F_i$'s, the above equations define K3 surfaces and $\sigma_{1,1}$, $\sigma_{4,1}$ have fixed locus of type $(1,1)$ and $(4,1)$ respectively.
\end{pro}
\proof
The fixed locus of $\sigma_{1,1}$ is the union of the curve $C_{1,1}$ of genus $3$ and an isolated point. By Riemann Hurwitz formula the intersection of $C_{1,1}$ with any elliptic curve on $X_{1,1}$ is either zero or three. Hence, by \cite[Theorem 5.2]{SD}, $\phi_{1,1}$ is an embedding of $X_{1,1}$ in $\Ps 3$.
Since $Im(\phi_{1,1})$ is invariant for the action of $\tilde \sigma_{1,1}$ we find that its equation is of the above form.
%Note that the fixed locus of $\tilde \sigma$ on $X$ is the union of $C=X\cap \pi$ and the point $p=(0,0,0,1)$.\\

According to Proposition \ref{transc}, the Picard lattice of $X_{4,1}$ is isomorphic to $U(3)\oplus A_2^4$.  Let $e,f$ be the basis of $U(3)$ as before and $e_i, e'_i$, $i=1,2,3,4$ the standard basis of $A_2^4$.  
The class $h=2e+f-\sum_{i=1}^4 (e_i+e'_i)$ gives an embedding as a smooth quartic surface in $\Ps 3$. The equations of $X_{4,1}$ and $\sigma_{4,1}$ can be determined as before.

The last statement follows from easy computations.
\qed
\begin{rem}
If $g(n,k)=3$ and $n\geq 2$ then $C_{n,k}$ is hyperelliptic by Proposition \ref{hyp},
 hence by \cite{SD} $\phi_{n,k}$ is a rational map of degree $2$ onto a quadric $Q$. Since $Q$ is invariant for $\tilde\sigma_{n,k}$ in Lemma \ref{proj}, then $Q$ is a cone and $\phi_{n,k}$ is branched along a quartic section $B$ and the vertex.
 The rational curves orthogonal to $C_{n,k}$ are contracted to the vertex of the cone and give a singular point of $B$ of type $A_1$ for $n=2$, $A_4$ for $n=3$ and $A_7$ for $n=4$. 
\end{rem}
\noindent The following was proved in \cite{D}.
\begin{pro}\label{2,1}
The surface $X_{2,1}$ is isomorphic to the double cover of $\Ps 2$ branched along a smooth plane sextic with equation
$$ F_6(x_0,x_1)+F_3(x_0,x_1)x_2^3+bx_2^6=0$$
and 
$$\sigma_{2,1}(x_0,x_1,x_2)=(x_0,x_1,\zeta x_2).$$
Conversely, for generic $F_i$'s, the above equation defines a K3 surface and $\sigma_{2,1}$ has fixed locus of type $(2,1)$.
\end{pro}
\proof
The morphism $\phi_{2,1}$ is a double cover of $\Ps 2$ branched along a smooth plane sextic by \cite{SD}.  Note that $\tilde \sigma_{2,1}$   (Lemma \ref{proj}) fixes pointwisely the line $x_2=0$ and the point $p=(0,0,1)$.

The last statement follows from easy computations.
\qed
\begin{rem}
If $g(n,k)=2$ and $n>2$ then it can be proved that $\phi_{n,k}$ is a degree two morphism to $\Ps 2$ branched along a plane sextic as in Proposition \ref{2,1} with $b=0$ i.e. $p=(0,0,1)$ is singular.
One can show that if $F_3$ has no multiple roots, then $p$ is of type $D_4$. If $F_3$ has exactly a double root $q$, then $p$ is of type 
$D_7$ if  $F_6(q)\not=0$ and of type $D_{10}$ if $F_6(q)=0$ is a simple root. This completes \cite[Proposition 3.3.18]{D}.
\end{rem}

\section{The moduli space}
We denote by $\mathcal M_{n,k}$ the moduli space of $K3$ surfaces with an order three non-symplectic automorphism with $n$ fixed points and $k$ fixed curves i.e. the space of  all pairs $(X_{n,k},\sigma_{n,k})$.
\begin{pro}\label{uniq} $\mathcal M_{n,k}$ is irreducible.
\end{pro}
\proof In propositions \ref{k>1}, \ref{ic}, \ref{qu} and \ref{2,1} we proved that the surface $X_{n,k}$ is isomorphic to the general element of an irreducible family of surfaces and that the automorphism $\sigma_{n,k}$ is uniquely determined. This gives the result.\qed\\

Let $\rho$ be an order $3$ isometry such that $(T,\rho)$ is a $\mathcal E^*$-lattice and consider the period domain $B_{\rho}$ as in section 3.
Since $T$ has signature $(2,2m-2)$ it is easy to see that $B_{\rho}$ is isomorphic to a $(m-1)$-dimensional complex ball.
Let $$\Gamma_{\rho}=\{\gamma \in O(T): \gamma\circ \rho=\rho \circ \gamma\}.$$ 
\begin{theorem}[\S 11, \cite{DK}] \label{mod} The generic point of
$\mathcal M_{\rho}=B_{\rho}/\Gamma_{\rho}$
parametrizes pairs $(X,\sigma)$ where $X$ is a K3 surface and $\sigma$ an order $3$ non-symplectic automorphism on $X$ with $\sigma^*=\rho$ up to isometries.
\end{theorem}

\begin{cor}
For any $n,k$ there is a unique isometry $\rho_{n,k}$ such that $T(n,k)$ is a $\mathcal E^*$-lattice and $\mathcal M_{n,k}$ is birational to $\mathcal M_{\rho_{n,k}}$.
\end{cor}
\proof Observe that $\mathcal M_{\rho}$ is irreducible. Then the result follows from Proposition \ref{uniq} and Theorem \ref{mod}. 
\qed\\

According to Theorem \ref{clas}, $T(3,0)=U\oplus U(3)\oplus A_2^{5}$. By Theorem \ref{RS} we also have an isomorphism $T(3,0)\cong A_2(-1)\oplus K_{12}$, where $K_{12}$ is the Coxeter-Todd lattice.

\begin{pro}\label{cox}
If $K_{12}$ is a primitive $\mathcal E$-sublattice of $T(n,k)$ then $(n,k)=(0,3)$.
\end{pro}
\proof
Let $B$ be the orthogonal complement of $K_{12}$ in $T=T(n,k)$.
Since $K_{12}$ is a unimodular $\mathcal E$-lattice then $K_{12}\oplus  B=T$. In fact, assume on the contrary that there exist $a\in K_{12}$, $b\in B$, $n\in \Z$ such that $(a+b)/n \in T/(K_{12}\oplus B)$.
Let $H$ be the hermitian form associated to a $\mathcal E$-lattice in Remark \ref{herm}, then the homomorphism
$$\phi:K_{12}\lra \mathcal E,\ x\mapsto H(x,(a+b)/n)=H(x,a/n)$$
belongs to $Hom(K_{12},\mathcal E)$ and not to $K_{12}$.

Since the lattice $K_{12}$ is negative definite and  $T$ has signature $(2,20-s)$ then $\rank B\geq 2$. 
Moreover $a(T)\geq 6$, hence $\rank T^{\perp}\geq 6$ and  $\rank B\leq 4$.

If $\rank B=2$ then $B$ is a rank one positive definite $\mathcal E^*$-lattice. A direct computation shows that $B\cong A_2(-1)$, hence $T\cong T(3,0)$. 
 
If $\rank B=4$ then $a(T)=6$ (since $\rank T^{\perp}=6$), hence $B$ is a unimodular lattice of signature $(2,2)$. By \cite[Corollary 1.13.3]{N2} $B\cong U\oplus U$.
 Assume that  $K_{12}\oplus U\oplus U$ admits a primitive embedding in $L_{K3}$. Then its orthogonal complement is a hyperbolic $3$-elementary lattice of rank $6$ with $a=6$. This lattice does not exist by Theorem \ref{RS}, hence this case can not occur.
 \qed

\begin{rem}
The Coxeter-Todd lattice also appears in connection to symplectic automorphisms of order $3$ on K3 surfaces. In \cite[Theorem 5.1]{GS} it is proved that the orthogonal complement of the fixed lattice of such automorphisms is isomorphic to $K_{12}$.
\end{rem}

\begin{theorem}\label{moduli}
The moduli space of K3 surfaces with a non-symplectic automorphism of order $3$ has three irreducible components which are the closures of
$$\mathcal M_{0,1},\ \mathcal M_{0,2},\ \mathcal M_{3,0}.$$
 \end{theorem}
\proof
By Proposition \ref{k>1} it is clear that all moduli spaces $\mathcal M_{n,k}$ with $k>1$ are in the closure of $\mathcal M_{0,2}$.

By Remark \ref{g4} a K3 surface in $\mathcal M_{0,1}$ is the triple cover of a quadric surface $Q$ branched along a genus $4$ curve $C$.
It is easy to check that, if $C$ has singular points of type $A_1$ then the triple cover of $Q$ branched along $C$ has rational double points of type $A_2$ and its minimal resolution  is still a K3 surface with an order $3$ non-symplectic automorphism.   
For example, if $C$ is generic with a node, then the associated K3 surface has an order three automorphism which fixes the proper transform of $C$ and one point (coming from the contraction of a component of the exceptional divisor over the node). In fact the Picard lattice isomorphic to $U(3)\oplus A_2$, where $A_2$ is generated by components of the exceptional divisor. Hence we are in case $n=1$, $k=1$.
Similarly, by taking $C$ with $n=1,\dots,4$ singularities of type $A_1$, we obtain cases $(n,1)$.

As a consequence of Proposition \ref{cox}, $\mathcal M_{3,0}$ is not contained in any other $\mathcal M_{n,k}$, hence it gives a maximal irreducible component of the moduli space.
\qed

\begin{rem}
In \S 4, \cite{K} S. Kond\=o describes the relation between the two components $\mathcal M_{0,1}$ and $\mathcal M_{0,2}$ (in fact the second component contains jacobians of the first one) and relates the groups $\Gamma_{0,1}, \Gamma_{0,2}$ with those appearing in Deligne-Mostow's list \cite{DM}.
\end{rem}

\bibliographystyle{plain}

\begin{thebibliography}{10}
 \bibitem{AN} V. Alexeev, V.V. Nikulin.
 \newblock Del Pezzo and K3 surfaces.
 \newblock MSJ Memoirs, Vol. 15, 2006.

\bibitem{AS} M.F. Atiyah, I.M. Singer.
\newblock {\em The Index of Elliptic Operators: III.}
\newblock Annals of Mathematics, 2nd Ser., Vol. 87, No. 3 (1968) 546--604.

\bibitem{Bu} D. A. Buell.
\newblock{ Binary quadratic forms: classical theory and modern computations}.
\newblock Springer-Verlag, 1989.


\bibitem{Br} G.E. Bredon.
\newblock {Introduction to compact Transformation Groups}.
\newblock Academic Press, New York, 1972.

\bibitem{Ca} H. Cartan.
\newblock {\em Sur les groupes de transformation analytiques.}
\newblock Act. Sci. Industrielles, 198, 1935.

\bibitem{CS} J.H. Conway, N.J.A. Sloane. 
\newblock Sphere Packing, Lattices and Groups. 
\newblock Springer-Verlag, New York, 1998.

\bibitem{CS1} J.H. Conway, N.J.A.  Sloane.
\newblock{\em The Coxeter-Todd lattice, the Mitchell group, and related sphere packings}.
\newblock Math. Proc. Camb. Phil. Soc., 93 (1983) 421--440.


\bibitem{D} J. Dillies.
\newblock{\em Automorphisms and Calabi-Yau Threefolds}.
\newblock{PhD thesis, University of Pennsylvania, 2006.}

\bibitem{DK} I. V. Dolgachev, S. Kond\=o. 
\newblock Moduli spaces of K3 surfaces and complex ball quotients.   \newblock {\em  Arithmetic and Geometry Around Hypergeometric Functions Lecture Notes of a CIMPA Summer School held at 
Galatasaray University, Istanbul, 2005}, Birkhauser Verlag Basel, Series: Progress in Mathematics, Vol. 260, 2007, 43--100. 
 
\bibitem{DM} P. Deligne, G.D. Mostow.
\newblock {\em Monodromy of hypergeometric functions and non-lattice integral monodromy}.
\newblock Publ. Math. IHES, 63 (1986), 5--89.

%\bibitem{DKV} I. Dolgachev, S. Kond\=o, B. van Geemen.
%\newblock{\em A complex ball uniformization of the moduli space of cubic surfaces via periods of {$K3$} surfaces},
%\newblock J. Reine Angew. Math. 588 (2005) 99--148.


\bibitem{F} W. Feit.
\newblock {\em Some lattices over {${\bf Q}(\surd -3)$}}.   
\newblock J. Algebra, 52, n.1 (1978) 248--263.

%\bibitem{Fr} Frame...

\bibitem{GS} A. Garbagnati, A. Sarti.
\newblock{\em Symplectic automorphisms of prime order on K3 surfaces}.
\newblock J. Algebra, 318 (2007) 323--350.


\bibitem{Ka} V.M. Kharlamov.
\newblock{\em The topological type of nonsingular surfaces in {$\mathbb P^3\mathbb R$}  of degree four},  
\newblock Funct. Anal. Appl.,10:4 (1976) 295Ð-304.  

\bibitem{K} S. Kond\=o.
\newblock{\em The moduli space of curves of genus $4$ and Deligne-Mostow's complex reflection groups}.
\newblock Advanced Studies in Pure Math., 36, Algebraic Geometry 2000, Azumino, 383--400.

\bibitem{Ko} S. Kond\=o.
\newblock{\em Automorphisms of algebraic K3 surfaces which act trivially on Picard groups}.
\newblock J. Math. Soc. Japan, 44, n.1 (1992).
 
\bibitem{Ku} V. S. Kulikov.
\newblock{\em Degenerations of K3 surfaces and Enriques surfaces}.
\newblock Izv. Akad. Nauk SSSR Ser. Mat., 41 (1977), 1008--1042.

\bibitem{M} R. Miranda.
\newblock{\em The basic theory of elliptic surfaces}.
\newblock ETS Editrice Pisa (1989).


 \bibitem{MO}
 N. ~Machida and K. ~Oguiso.
 \newblock {\em On {$K3$} surfaces admitting finite non-symplectic group
              actions}.
              \newblock {\em J. Math. Sci. Univ. Tokyo,} vol. 5, n.2 (1998) 273--297.



\bibitem{Na} Y. Namikawa.
\newblock{\em Periods of Enriques surfaces},
\newblock Math. Ann., 270 (1985), 201--222.

\bibitem{N1}
V.V. Nikulin.
\newblock{\em Finite groups of automorphisms of K\"{a}hlerian surfaces of type $K3$}.
\newblock Moscow Math. Sod., 38 (1980), 71--137.

\bibitem{N2}
V.V. Nikulin
\newblock{\em Integral symmetric bilinear forms and some of their applications}.
\newblock Math. USSR Izv., 14 (1980), 103--167.

\bibitem{N3}
V.V. Nikulin
\newblock{\em Factor groups of groups of the automorphisms of hyperbolic forms with respect to subgroups generated by 2-reflections}.
\newblock J. Soviet Math., 22 (1983), 1401--1475.


\bibitem{OZ2}
K. Oguiso and D-Q. Zhang.
\newblock{\em On Vorontsov's theorem on K3 surfaces with non-symplectic group actions}.
\newblock Proc. Amer. Math. Soc. 128 (2000), n.6, 1571--1580.

\bibitem{OZ1}
K. Oguiso and D-Q. Zhang.
\newblock{\em K3 surfaces with order 11 automorphisms}.
\newblock math.AG/9907020.


\bibitem{PP}
U. Persson and H. Pinkham.
\newblock{\em Degenerations of surfaces with trivial canonical bundle}.
\newblock Ann. Math., 113 (1981), 45--66.

 
\bibitem{RS}
A. N. Rudakov and I. Shafarevich.
\newblock{\em Surfaces of type K3 over fields of finite characteristic.}
\newblock In: I. Shafarevich, Collected mathematical papers, Springer, Berlin (1989), 657--714.

\bibitem{SD}
B.~Saint-Donat.
\newblock {\em Projective models of K3 surfaces.}
\newblock American Journal of Math. 96(4) (1974), 602--639.

\bibitem{S}
H. Sterk.
\newblock{\em Finiteness results for algebraic K3 surfaces}.
\newblock Mathematische Zeitschrift, 189 (1985), 507--513.


\bibitem{Z1} D-Q. Zhang.
\newblock{\em Automorphisms of K3 surfaces}.
\newblock AMS/IP Studies in advanced Mathematics, 39 (2007), 379--392.

\bibitem{Z2} D-Q. Zhang.
\newblock{\em Quotients of K3 surfaces modulo involutions}.
\newblock Japanese J. Mathematics, to appear (1999).

\bibitem{Z3} D-Q. Zhang.
\newblock {\em Normal algebraic surfaces  with trivial tricanonical divisors.}
\newblock Publ. RIMS, Kyoto Univ., 33 (1997), 427--442.



\end{thebibliography}

\end{document}